\documentclass{amsart}
\usepackage{amscd}
\usepackage{amssymb}
\usepackage{amsmath}

\usepackage[all]{xy}
\theoremstyle{plain}
\newtheorem{theorem}{Theorem}[section]
\newtheorem{corollary}[theorem]{Corollary}

\newtheorem{lemma}[theorem]{Lemma}

\newtheorem{proposition}[theorem]{Proposition}
\newtheorem{question}[theorem]{Question}

\theoremstyle{definition}
\newtheorem{definition}[theorem]{Definition}
\newtheorem{example}[theorem]{Example}

\newtheorem{remark}[theorem]{Remark}
\newtheorem{Definitions and Notation}[theorem]{Definitions and
Notation}

\numberwithin{equation}{section}

\newcommand{\f}[1]{\ensuremath{\mathfrak{#1}}}
\newcommand{\locoho}[3]{\ensuremath{\mathrm{H}_{#1}^{#2}\left(#3\right)}}

\newcommand{\sse}{\subseteq}

\newcommand{\soc}[1]{\ensuremath{\mathrm{S}\left(#1\right)}}
\newcommand{\socdim}[1]{\ensuremath{\f{S}\left(#1\right)}}
\renewcommand{\index}[3]{\ensuremath{\mathrm{N}_{#1}\left(#2;#3\right)}}
\renewcommand{\hom}[3]{\ensuremath{\mathrm{Hom}_{#1}\left(#2,#3\right)}}
\renewcommand{\colon}[3]{\ensuremath{\left(#1:_{#2} #3\right)}}
\newcommand{\ext}[4]{\ensuremath{\mathrm{Ext}_{#1}^{#2}\left(#3,#4\right)}}

\newcommand{\type}[2]{\ensuremath{\mathrm{r}_{#1}\left(#2\right)}}
\newcommand{\classicaltype}[2]{ \ensuremath{\overline{r}_{#1} \left(#2\right) }}
\newcommand{\idealize}[2]{\ensuremath{#1 \ltimes #2}}
\newcommand{\ann}[2]{\ensuremath{\mathrm{ann}_{#1}\left(#2\right)}}
\newcommand{\tensor}[3]{\ensuremath{#1 \otimes_{#2} #3}}
\newcommand{\length}[2]{\ensuremath{\lambda_{#1} \left(#2\right)}}
\newcommand{\minimal}[2]{\ensuremath{\mu_{#1}(#2)}}
\newcommand{\complete}[1]{\ensuremath{\hat{#1}}}
\newcommand{\reduction}[2]{\ensuremath{\mathrm{red}_{#1} \left(#2\right)}}
\newcommand{\globalreduction}[1]{\ensuremath{\mathrm{red} \left(#1\right)}}
\newcommand{\unmixed}[2]{\ensuremath{\text{U}\left(#1; #2 \right)}}

\DeclareMathOperator{\depth}{depth}
\DeclareMathOperator{\Dim}{dim}

\begin{document}

\title[The index of reducibility of parameter ideals]
{The index of reducibility of parameter ideals\\
in low dimension}

\author{Mark W. Rogers}

\address{Mathematics Department,
Purdue University, West Lafayette, IN 47907-2066}

\email{mrogers@math.purdue.edu}


\begin{abstract}
In this paper we present results concerning the following
question: If $M$ is a finitely-generated module with finite local
cohomologies over a Noetherian local ring $(A, \f{m})$, does there
exist an integer $\ell$ such that every parameter ideal for $M$
contained in $\f{m}^{\ell}$ has the same index of reducibility? We
show that the answer is \emph{yes} if $\Dim M = 1$ or if $\Dim M =
2$ and $\depth M > 0$.  This research is closely related to work
of Goto-Suzuki and Goto-Sakurai; Goto-Sakurai have supplied an
answer of \emph{yes} in case $M$ is Buchsbaum.
\end{abstract}
\maketitle

\section{Introduction}

Let $A$ be a $d$-dimensional Noetherian local ring with maximal
ideal \f{m} and residue field $k=A/\f{m}$, and let $M$ be a
finitely generated $A$-module. Recall that a submodule of $M$ is
called \emph{irreducible} if it cannot be written as the
intersection of two larger submodules.  It is well known that
every submodule $N$ of $M$ can be expressed as an irredundant
intersection of irreducible submodules, and that the number of
irreducible submodules appearing in such an expression depends
only on $N$ and not on the expression
\cite[p.~92-3]{Sharpe.Vamos}.

For an ideal $I$ of $A$, we say that $I$ is \emph{cofinite} on $M$
if the module $M/IM$ has finite length.  For an ideal $I$ of $A$
which is cofinite on $M$, the \emph{index of reducibility} of $I$
on $M$ is defined as the number of submodules appearing in an
irredundant expression of $IM$ as an intersection of irreducible
submodules of $M$.  We denote the index of reducibility of $I$ on
$M$ by \index{A}{I}{M}.

The smallest number of generators of an ideal which is cofinite on
$M$ is the dimension of $M$, and a cofinite ideal having this
minimal number of generators is called a \emph{parameter ideal}
for $M$.  Our interest in the index of reducibility of parameter
ideals stems from the relationship with the Cohen-Macaulay and
Gorenstein properties.  In 1956, D. G. Northcott proved that in a
Cohen-Macaulay local ring, the index of reducibility of any
parameter ideal depends only on the ring \cite[Theorem 3]{Northc}.
This result extends to modules, and the common index of
reducibility of parameter ideals for a Cohen-Macaulay module $M$
is called the \emph{(Cohen-Macaulay)} \emph{type} of $M$.  We
denote the type of $M$ by \type{A}{M}.

As a partial converse, Northcott along with D. Rees proved in
\cite{NR-princ} that if every parameter ideal of $A$ is
irreducible, then $A$ is Cohen-Macaulay.  This provides an
attractive characterization of a Gorenstein local ring as a local
ring in which every parameter ideal is irreducible.

One might suspect that the Cohen-Macaulay property of a local ring
is characterized by the constant index of reducibility of
parameter ideals.  However, in 1964, S. Endo and M. Narita in
\cite{Endo.Narita} gave an example of a Noetherian local ring in
which the index of reducibility of each parameter ideal is two,
and yet the ring is not Cohen-Macaulay.

In 1984, S. Goto and N. Suzuki generalized the example of
Endo-Narita, as well as undertaking a study of the supremum of the
index of reducibility of parameter ideals for $M$
\cite{Goto.Suzuki}. We refer to this supremum as the
\emph{Goto-Suzuki} type (\emph{GS-type}) of $M$, and denote it by
\classicaltype{A}{M}. In the case where $M$ is Cohen-Macaulay, the
GS-type \classicaltype{A}{M} is equal to the type \type{A}{M}.
However, Goto-Suzuki provide examples where the GS-type of a
Noetherian local ring is infinity.

We introduce some terminology in order to state one of the main
results of Goto-Suzuki.  We denote the $i$th local cohomology
module of $M$ with respect to \f{m} by \locoho{\f{m}}{i}{M}. We
say that $M$ has \emph{finite local cohomologies} if the modules
\locoho{\f{m}}{i}{M} have finite length for $i\neq d$. We use
\length{A}{M} to denote the length of $M$, and we set
$\socdim{M}=\length{A}{\hom{A}{k}{M}}$, the \emph{socle dimension}
of $M$. We let \minimal{A}{M} denote the minimal number of
generators of $M$, and let $E$ denote the injective hull of the
residue field $k$.  The main result of Goto-Suzuki concerning the
GS-type of a module having finite local cohomologies is the
following:

\begin{theorem}[Goto-Suzuki]
Let $M$ be a finitely generated $d$-dimensional $A$-module with
finite local cohomologies.  Then we have the following
inequalities:
\begin{equation}\label{eq:gotosuzuki}
\sum_{i=0}^{d} \binom{d}{i} \socdim{\locoho{\f{m}}{i}{M}} \leq
\classicaltype{A}{M} \leq \sum_{i=0}^{d-1} \binom{d}{i}
\length{A}{\locoho{\f{m}}{i}{M}} + \minimal{\complete{A}}{K},
\end{equation}
where \complete{A} is the \f{m}-adic completion of $A$ and $K =
\hom{A}{\locoho{\f{m}}{d}{M}}{E}$ is the canonical module of the
completion of $M$.
\end{theorem}

\begin{proof}
See \cite[Theorems~2.1 and 2.3]{Goto.Suzuki}.
\end{proof}

In 1994, Kawasaki used this result to determine conditions under
which a module having finite GS-type is Cohen-Macaulay. The main
result of Kawasaki in \cite{Kawasaki} is that if $A$ is the
homomorphic image of a Cohen-Macaulay ring, then $M$ is
Cohen-Macaulay if and only if \classicaltype{A}{M} is finite,
$M_{\f{p}}$ is Cohen-Macaulay for all primes in the support of $M$
with $\dim M_{\f{p}} < \classicaltype{A}{M}$, and all the
associated primes of $M$ have the same dimension.

In two recent papers \cite{Goto.SakuraiI, Goto.SakuraiII}, Goto
with H. Sakurai has returned to the study of the index of
reducibility of parameter ideals in order to investigate when the
equality $I^2=QI$ holds for a parameter ideal $Q$ in $A$, where
$I=\colon{Q}{A}{\f{m}}$. According to earlier research of A.
Corso, C. Huneke, C. Polini, and W. Vasconcelos
\cite{Corso.Huneke.Vasconcelos, Corso.Polini,
Corso.Polini.Vasconcelos}, this equality holds for all parameter
ideals $Q$ in case $A$ is a Cohen-Macaulay ring which is not
regular. Goto-Sakurai generalize this and say that if $A$ is a
Buchsbaum ring whose multiplicity is greater than 1, then the
equality $I^2 = QI$ holds for any parameter ideal whose index of
reducibility is the GS-type of $A$.  Thus, if $A$ is a Buchsbaum
ring whose multiplicity is greater than 1 and if $A$ has constant
index of reducibility of parameter ideals, then the equality
$I^2=QI$ holds for all parameter ideals $Q$.

Most pertinent to our discussion is Corollary 3.13 of Goto-Sakurai
in \cite{Goto.SakuraiI}, which states that if $A$ is a Buchsbaum
ring of positive dimension, then there is an integer $\ell$ such
that the index of reducibility of $Q$ is independent of $Q$ and
equals \classicaltype{A}{A} for all parameter ideals $Q \sse
\f{m}^{\ell}$.  In view of this, it is natural to ask the
following question:

\begin{question}
Suppose $(A, \f{m})$ is a Noetherian local ring having finite
local cohomologies.  Is there an integer $\ell$ such that the
index of reducibility of any parameter ideal contained in
$\f{m}^{\ell}$ is the same?
\end{question}

We note that if $A$ is a ring satisfying the hypothesis of the
question, and if the answer to the question is \emph{yes} for $A$,
then the common index of reducibility of parameter ideals in high
powers of the maximal ideal is equal to the lower bound of
Goto-Suzuki: $\sum_{i = 0}^{d} \binom{d}{i}
\socdim{\locoho{\f{m}}{i}{A}}$. This is because implicit in the
proof of the lower bound \cite[Theorem 2.3]{Goto.Suzuki} we find
that given any system of parameters $x_1$, \ldots, $x_d$ for $A$,
there are integers $n_i$ such that the parameter ideal
$(x_{1}^{n_1}, \ldots, x_{d}^{n_d})A$ has index of reducibility
$\sum_{i = 0}^{d} \binom{d}{i} \socdim{\locoho{\f{m}}{i}{A}}$.

The main result of the current paper is the following:

\begin{theorem}
Let $(A, \f{m})$ be a Noetherian local ring and let $M$ be a
finitely-generated $A$-module of dimension $d \leq 2$.  Suppose
either $M$ has dimension 1, or $M$ has finite  local cohomologies
and depth at least one.  Then there exists an integer $\ell$ such
that for every parameter ideal \f{q} for $M$ contained in
$\f{m}^{\ell}$, the index of reducibility of \f{q} on $M$ is
independent of \f{q} and is given by
\begin{equation}
\index{A}{\f{q}}{M} = \sum_{i = 0}^d \binom{d}{i}
\socdim{\locoho{\f{m}}{i}{M}}.
\end{equation}
\end{theorem}

\begin{proof}
See Theorem \ref{MainTheorem} and \ref{MainTheorem2}.
\end{proof}

As a corollary of this result we prove that a Noetherian local
ring $A$ of dimension at most 2 is Gorenstein if and only if $A$
has finite local cohomologies, and inside every power of the
maximal ideal of $A$ there exists an irreducible parameter ideal.
A Noetherian local ring having the property that every power of
its maximal ideal contains an irreducible parameter ideal is
called an \emph{approximately Gorenstein} ring, or is said to have
\emph{small cofinite irreducibles (SCI)}.  The author's original
motivation for studying the index of reducibility of parameter
ideals comes from questions that arose while studying M.
Hochster's paper \cite{Hochster} exposing the relationship between
modules having SCI and the condition that cyclic purity implies
purity.

We present an example of a complete Noetherian local ring of
dimension $d$ and depth $d - 1$ ($d > 1$), such that
\locoho{\f{m}}{d - 1}{A} is not finitely generated, and such that
in every power of the maximal ideal there is a parameter ideal
with index of reducibility 2 and a parameter ideal with index of
reducibility 3. This example is obtained as an idealization, so we
present several basic results relating the index of reducibility
with an idealization. Using a result of C. Lech \cite{Lech}, we
are able to obtain such an example among Noetherian local domains.

The current paper is part of the author's PhD dissertation under
the guidance of Professor W. Heinzer at Purdue University.  It is
the author's pleasure to thank Professor Heinzer for all of his
guidance. The author is also pleased to thank the referee for very
insightful and generous remarks concerning Section
\ref{Dimension2}.

\section{Background and Dimension One}

We begin with some terminology.

\begin{Definitions and Notation}  Let $A$ denote a Noetherian local ring with
maximal ideal \f{m} and residue field $k = A/\f{m}$, let $I$
denote an ideal of $A$, and let $M$ denote a finitely generated
$A$-module.
\begin{enumerate}
\item The \emph{socle} of $M$ is defined to be $\soc{M} =
\colon{0}{M}{\f{m}}$.  Note that $\soc{M} \cong \hom{A}{k}{M}$.
The socle is naturally a vector space over $k$, and we denote its
dimension by \socdim{M}.

\item Suppose $M$ has depth $t$.  The \emph{type} of $M$ is
defined by
\begin{equation*}
\type{A}{M} = \Dim_k \ext{A}{t}{k}{M}.
\end{equation*}
Note that if $I$ is cofinite on $M$, then $\index{A}{I}{M} =
\type{A}{M/IM}$, and if $M$ has depth 0, then $\socdim{M} =
\type{A}{M}$.  See \cite[p.~13]{Bruns.Herzog} for more information
on the type of a module.

\item Given an $A$-algebra $B$, we define $M_B =
\tensor{M}{A}{B}$. Then $(A/I)_B \cong B/IB$, and if $B$ is flat
over $A$, then for any submodule $N$ of $M$, we have $M_B /N_B
\cong (M/N)_B$.

\item We say that an ideal $J$ of $A$ is a \emph{reduction} of $I$
if $J \sse I$ and there is some integer $n$ such that $I^{n+1} = J
I^n$.  We say that $J$ is a \emph{minimal reduction} of $I$ if
there is no reduction of $I$ properly contained in $J$.  The
\emph{reduction number} of $I$ with respect to $J$ is defined as
\begin{equation*}
\reduction{J}{I}=\min \{n :\, I^{n+1} = JI^n \}.
\end{equation*}
We define the \emph{reduction number} of $I$ to be
\begin{equation*}
\globalreduction{I} = \min \{ \reduction{J}{I} :\, J \text{ is a
minimal reduction of }I \}.
\end{equation*}
\end{enumerate}
\end{Definitions and Notation}

We must reduce to the case that the residue field is infinite in
order to obtain a principal reduction of the maximal ideal. Thus
our first result is a basic lemma concerning the behavior of the
index of reducibility under a flat, local change of base.  We omit
the proof of Lemma \ref{flat.base.change}, but cite
\cite[Theorem~19.1]{Nagata},
\cite[Proposition~1.2.16]{Bruns.Herzog} and
\cite[Theorem~7.4]{Matsumura} for relevant related results.

\begin{lemma}\label{flat.base.change}
Suppose $A$ and $B$ are Noetherian local rings with maximal ideals
\f{m} and \f{n}, respectively, let $M$ be a finitely generated
$A$-module, and suppose $B$ is a flat $A$-algebra for which $\f{n}
= \f{m}B$.  Let $I$ be an ideal of $A$ and let $N$ be a submodule
of $M$.
\begin{enumerate}
\item If $M$ has finite length, then $M_B$ has finite length, and
$\length{B}{M_B} = \length{A}{M}$.\label{finite.MB}

\item If $I$ is an ideal of $A$ which is cofinite on $M$, then
$IB$ is cofinite on $M_B$.

\item \label{flat.index} If $I$ is an ideal of $A$ which is
cofinite on $M$, then
\begin{equation*}
\index{B}{IB}{M_B} = \index{A}{I}{M}.
\end{equation*}

\item $\colon{N}{M}{I}_B \cong \colon{N_B}{M_B}{IB}.$
\label{flat.colon}

\item \label{flat.locoho} $\locoho{\f{n}}{0}{M_B} =
\locoho{\f{m}}{0}{M}_B .$
\end{enumerate}
\end{lemma}

Now we are prepared to prove that parameter ideals in high powers
of a one dimensional Noetherian local ring all have the same index
of reducibility.

\begin{theorem}\label{MainTheorem}
Suppose $A$ is a Noetherian local ring with maximal ideal \f{m},
and let $M$ be a finitely generated $A$-module of dimension 1. Set
$W=\locoho{\f{m}}{0}{M}$. There is an integer $\ell$ such that the
index of reducibility of any parameter ideal $\f{q} \sse
\f{m}^{\ell}$ is independent of \f{q} and is given by
\[
\index{A}{\f{q}}{M} = \socdim{M} + \type{A}{M/W}.
\]

When $M=A$ and $k$ is infinite, the integer $\ell$ may be taken to
be
\begin{equation*}
\ell = \max \{c, d \}+1
\end{equation*}
where $c = \globalreduction{\f{m}}$ and $d$ is the smallest
integer with $\f{m}^d \cap W = 0$.
\end{theorem}

\begin{proof}
First we show that it suffices to prove the theorem in the case
where the residue field $k$ is infinite. Suppose $B$ is any local
flat $A$-algebra whose maximal ideal \f{n} is extended from $A$
and whose residue field is infinite (for instance,
$B=A[x]_{\f{m}A[x]}$ for an indeterminate $x$). Then according to
Theorem~A.11 on p.~415 of \cite{Bruns.Herzog}, the $B$-module
$M_B$ has dimension 1. Therefore we may apply the theorem to
$M_B$; let $\ell$ be the integer guaranteed by the theorem. Set $V
= \locoho{\f{n}}{0}{M_B}$; then according to
Part~\ref{flat.locoho} of Lemma~\ref{flat.base.change}, we have $V
= W_B$.

Suppose \f{q} is a parameter ideal for $M$ that is contained in
$\f{m}^{\ell}$. Then $\f{q}B$ is a parameter ideal for $M_B$ which
is contained in $\f{n}^{\ell}$, so we have
\begin{equation}\label{infinite.case}
\index{B}{\f{q}B}{M_B} = \socdim{M_B} + \type{B}{M_B /V}.
\end{equation}

According to Part~\ref{flat.index} of
Lemma~\ref{flat.base.change}, the left side of
equality~\ref{infinite.case} equals \index{A}{\f{q}}{M}, so it
remains to see that $\socdim{M_B} = \socdim{M}$ and $\type{B}{M_B
/V} = \type{A}{M/W}$. For the first equality, if $W = 0$ then $V =
0$ and we have nothing to show.  Otherwise, $M$ has depth 0 and
thus according to Proposition~1.2.16 part~(a) from
\cite{Bruns.Herzog}, so does $M_B$, so by part~(b) of the same
Proposition, we have
\begin{equation*}
\socdim{M_B} = \type{B}{M_B} = \type{A}{M} = \socdim{M}.
\end{equation*}

To see the equality $\type{B}{M_B /V} = \type{A}{M/W}$, we first
note that
\begin{equation*}
M_B /V = M_B /W_B \cong (M/W)_B.
\end{equation*}
Hence, by another application of Proposition~1.2.16 part~(b) from
\cite{Bruns.Herzog}, we have
\begin{equation*}
\type{B}{M_B /V} = \type{B}{(M/W)_B} = \type{A}{M/W}.
\end{equation*}
Thus we see that we may assume the residue field $k$ is infinite.

Our second task is to see that we may replace the ring $A$ by the
ring $C=A/I$, where $I$ is the annihilator of $M$.  Let \f{q}
denote an ideal of $A$ which is a parameter ideal for $M$.  The
$A$-module $M$ becomes an $C$-module in a natural way, and the
submodule structure of $M$ remains unchanged.  The maximal ideal
of $A$ extends to that of $C$, the residue field of $C$ is $k$,
and $\f{q}C$ is an ideal of $C$ which is a parameter ideal for
$M$. If $N$ is any finitely generated $A$-module annihilated by
$I$, then as sets we have $\hom{A}{k}{N} = \hom{C}{k}{N}$, and we
easily check that this equality is actually an isomorphism of
$A$-modules. Thus \socdim{N} does not depend on whether we view
$N$ as an $A$-module or a $C$-module. Since the numbers
\index{A}{\f{q}}{M}, \socdim{M}, and \type{A}{M/W} are all
calculated as socle dimensions of quotients of $M$, they do not
change when we view $M$ as an $C$-module. Thus we replace $A$ by
$C$ and assume that $A$ has Krull dimension one.

If $W = 0$, then $M$ is Cohen-Macaulay, $\socdim{M} = 0$, and for
a parameter ideal $\f{q} = aA$ we have
\begin{equation}
\index{A}{\f{q}}{M} = \socdim{M/aM} = \Dim_k \ext{A}{1}{k}{M} =
\type{A}{M}.
\end{equation}
Thus the proof is complete in this case.

Now suppose $W \neq 0$.  Since $W$ has finite length and
\begin{equation*}
\bigcap_{i = 1}^{\infty} (\f{m}^n M \cap W) = 0,
\end{equation*}
there is some integer $d$ with $\f{m}^d M \cap W = 0$.

Since the residue field is infinite, \f{m} has a principal
reduction; i.e., there is an element $x \in \f{m}$ and an integer
$c \geq 1$ such that $\f{m}^{c+1} = x\f{m}^{c}$
(\cite[Corollary~4.6.10, p.~191]{Bruns.Herzog}). Set $\ell = \max
\{c,d\} + 1$. Then we have arrived at a situation where any
parameter $a$ for $M$ which is in $\f{m}^{\ell}$ is of the form $a
= xy$ with $y \in \f{m}^d$.  Furthermore, since we have assumed
$M$ is faithful, the parameters for $M$ are just the parameters
for $A$ \cite[Proposition~10.8, p.~237]{Eisenbud}.  Since the
parameters for the one-dimensional ring $A$ are those elements not
in any minimal prime ideal of $A$, we see that if $a = xy$ is a
parameter for $M$, then so are $x$ and $y$.

Since $aM \cap W = 0$, we have $(W + aM)/aM \cong W$, so, as in
the proof of \cite[Proposition (2.4)]{Goto.Suzuki}, we see that
the top row in the following commutative diagram is exact:
\begin{equation}
\xymatrix{ 0 \ar[r] & W \ar[r] & \frac{M}{aM} \ar[r]^{\beta} &
\frac{M}{aM + W}
 \ar[r] & 0\\
 & & & \frac{M}{xM + W} \ar[u]_{g} \ar[lu]^{f}   &
}
\end{equation}

The maps $f$ and $g$ are each given by multiplication by $y$. The
map $f$ is well-defined since $y(xM + W) = aM + yW = aM$, and $g$
is just the composition $g = \beta f$.

To see that $g$ is injective, view $g$ as multiplication by $y$ as
follows:
\[
\xymatrix{
 \frac{M/W}{x(M/W)} \ar[r]^{g = \cdot y} &\frac{M/W}{a(M/W)}}
\]
Since $M/W$ is a Cohen-Macaulay module and $y$ is a parameter for
$M$, and thus for $M/W$, we have that $y$ is regular on $M/W$.
Hence
\begin{equation*}
\colon{xy(M/W)}{M/W}{y} = x(M/W),
\end{equation*}
so that $g$ is injective.

Information on the dimension of the socles is obtained by applying
the functor \hom{A}{k}{-}; what results is the following exact
commutative diagram, where the maps induced by $f$ and $g$ are
still injective:
\begin{equation}\label{diagram}
\xymatrix{ 0 \ar[r] & \hom{A}{k}{W} \ar[r] &
\hom{A}{k}{\frac{M}{aM}} \ar[r]^{\beta^*} &
\hom{A}{k}{\frac{M}{aM+W}}
 \\
 & & & \hom{A}{k}{\frac{M}{xM+W}} \ar[u]_{\hom{A}{k}{g}} \ar[lu]^{\hom{A}{k}{f}}   &
}
\end{equation}
Here we use $\beta^*$ to denote \hom{A}{k}{\beta}.

The important point is that since $x$ and $a$ are still parameters
for the Cohen-Macaulay module $M/W$, we have
\begin{eqnarray*}
\Dim_k \hom{A}{k}{\frac{M}{aM+W}} &=& \Dim_k
\hom{A}{k}{\frac{M/W}{a(M/W)}} \\
&=& \index{A}{aA}{M/W} \\
&=& \type{A}{M/W},
\end{eqnarray*}
and similarly
\begin{equation*}
\Dim_k \hom{A}{k}{\frac{M}{aM+W}} = \type{A}{M/W}.
\end{equation*}
The map \hom{A}{k}{g} is an injection of $k$-vector spaces, each
of dimension \type{A}{M/W}, hence this map is an isomorphism. From
the surjectivity of this map, it follows that $\beta^*$ is
surjective, so that the top row of Diagram~\ref{diagram} is exact.
Since the length of the middle module is \index{A}{aA}{M} and the
left-hand module is the socle of $W$, which is the socle of $M$,
we complete the proof in dimension one using the additivity of
length.
\end{proof}

\section{Dimension 2}\label{Dimension2}

We thank the referee for pointing out the fascinating technique in
this section, and for pointing out the recent work of Goto-Sakurai
\cite{Goto.SakuraiI, Goto.SakuraiII}.  Goto-Sakurai successfully
apply this technique to the case of Buchsbaum local rings of
arbitrary dimension \cite[Corollary 3.13]{Goto.SakuraiI}.

Suppose $\f{q} = (x_1, \ldots, x_d)A$ is an ideal generated by a
system of parameters $x_1$, \ldots, $x_d$ for a finitely-generated
module $M$ of dimension $d$ over a Noetherian local ring $(A,
\f{m})$.  We may form a direct system of modules by setting $M_i =
M/(x_{1}^i, \ldots, x_{d}^i)M$ and defining maps from $M_i \to
M_{i+1}$ given by multiplication by $x_1 \cdot \cdots \cdot x_d$.
It is known that the direct limit of this system is
\locoho{\f{m}}{d}{M} \cite[Theorem 3.5.6]{Bruns.Herzog}. From this
point of view, we see that there is a canonical homomorphism from
$M/\f{q}M$ to \locoho{\f{m}}{d}{M}.

For an ideal $I$ of a Noetherian local ring $(A,\f{m})$ and a
finitely-generated module $M$ we define \unmixed{I}{M} to be the
unmixed component of the submodule $IM$ of $M$; that is,
\unmixed{I}{M} is the intersection of the primary components of
the submodule $IM$ whose associated primes have maximal dimension,
equal to $\Dim M/IM$.

In the course of our proof for dimension 2 we will need to mention
several generalizations of the notion of a regular sequence. These
definitions can be found in the appendix of \cite{Stuckrad.Vogel},
which is a good source for information concerning modules having
finite local cohomologies.

\begin{definition}
Suppose $(A, \f{m})$ is a Noetherian local ring, $M$ is a
finitely-generated $A$-module of dimension $d > 0$, and let \f{a}
be an \f{m}-primary ideal.  A system of elements $x_1$, \ldots,
$x_r$ is called an \emph{\f{a}-weak $M$-sequence} if
\begin{equation*}
\colon{(x_1, \ldots, x_{i - 1})M}{M}{x_i} \sse \colon{(x_1,
\ldots, x_{i - 1})M}{M}{\f{a}}
\end{equation*}
for all $i=1$, \ldots, $r$.

Let $x_1$, \ldots, $x_d$ be a system of parameters for $M$ and let
$\f{q}=(x_1, \ldots, x_d)A$.  We say that $x_1$, \ldots, $x_d$ is
a \emph{standard system of parameters} of $M$ if $x_{1}^{n_1}$,
\ldots, $x_{d}^{n_d}$ is a \f{q}-weak $M$-sequence for all $n_1$,
\ldots, $n_d \geq 1$.

We say that \f{a} is a \emph{standard ideal} with respect to $M$
if every system of parameters of $M$ contained in \f{a} is a
standard system of parameters of $M$.
\end{definition}

We begin by isolating a general statement concerning the index of
reducibility of parameter ideals in the case where $M$ has finite
local cohomologies.

\begin{proposition}\label{T:lower bound}
Let $(A, \f{m},k)$ be a Noetherian local ring and let $M$ be a
finitely-generated $d$-dimensional $A$-module with $d>0$ such that
$M$ has finite local cohomologies.  There exists an integer $\ell$
such that for every parameter ideal $\f{q} = (x_1, \ldots, x_d)$
of $M$, if $\f{q} \sse \f{m}^{\ell}$ then the index of
reducibility of \f{q} on $M$ is given by
\begin{equation}
\index{A}{\f{q}}{M} = \socdim{\sum_{i = 1}^{d} \frac{U_i + x_i
M}{\f{q}M}} + \socdim{\locoho{\f{m}}{d}{M}},
\end{equation}
where $U_i = \unmixed{(x_1, \ldots, \widehat{x_i}, \ldots,
x_d)A}{M}$.
\end{proposition}

\begin{proof}
Since $M$ has finite local cohomologies, we have a standard ideal
\f{a} for $M$ \cite[Corollary 18, p. 264]{Stuckrad.Vogel}; hence
every system of parameters of $M$ contained in \f{a} is an
\f{a}-weak $M$-sequence \cite[Theorem 20, p. 264]{Stuckrad.Vogel}.
Thus, given any system of parameters $x_1$, \ldots, $x_d$ of $M$
contained in \f{a} and any integer $n \geq 1$, we have by
\cite[Lemma 23]{Stuckrad.Vogel}
\begin{align}
& \colon{(x_{1}^{n+1}, \ldots, x_{d}^{n+1})M}{M}{(x_1 \cdot \cdots
\cdot x_d)^n} \\
& \qquad = (x_1, \ldots, x_d)M + \sum_{i = 1}^d \colon{(x_1,
\ldots, \widehat{x_i}, \ldots, x_d)M}{M}{\f{a}}.\notag
\end{align}
We note that this right hand side is equal to
\begin{equation}
\sum_{i = 1}^d \left( \colon{(x_1, \ldots, \widehat{x_i}, \ldots,
x_d)M}{M}{x_i} + x_i M \right).
\end{equation}
Set $U_i = \unmixed{(x_1, \ldots, \widehat{x_i}, \ldots,
x_d)A}{M}$.  Since $U_i = \colon{(x_1, \ldots, \widehat{x_i},
\ldots, x_d)M}{M}{x_i}$, we have that
\begin{equation}\label{kernel}
\colon{(x_{1}^{n+1}, \ldots, x_{d}^{n+1})M}{M}{(x_1 \cdot \cdots
\cdot x_d)^n} = \sum_{i = 1}^d (U_i + x_i M).
\end{equation}

An important component of this proof is that according to
\cite[Lemma 3.12]{Goto.SakuraiI}, we may choose $\ell$ large
enough so that for any ideal \f{q} generated by a system of
parameters for $M$, if $\f{q} \sse \f{m}^\ell$ then the canonical
map $\phi : M/\f{q}M \to \locoho{\f{m}}{d}{M}$ is surjective on
the socles; that is, \hom{A}{k}{\phi} is surjective.
 We also require $\ell$ to be large enough so that
 $\f{m}^{\ell} \sse \f{a}$.

Let $\f{q} = (x_1, \ldots, x_d)$ be a parameter ideal of $M$
contained in $\f{m}^\ell$ and set $U_i = \unmixed{(x_1, \ldots,
\widehat{x_i}, \ldots, x_d)A}{M}$. Let $K$ denote the kernel of
the canonical map $\phi$ from $M/\f{q}M$ to \locoho{\f{m}}{d}{M}.
According to the definition of the direct limit, we have
\begin{equation}
K = \frac{\cup_{n \geq 1} \colon{(x_{1}^{n+1}, \ldots,
x_{d}^{n+1})M}{M}{(x_1 \cdot \cdots \cdot x_d)^n}}{\f{q}M};
\end{equation}
thus by Equation \ref{kernel} we see that
\begin{equation}
K = \sum_{i = 1}^d \frac{U_i + x_i M}{\f{q}M}.
\end{equation}

Now all that is left is to apply the socle functor \hom{A}{k}{-}
to the exact sequence
\begin{equation}
0 \to K \to M/\f{q}M \to \locoho{\f{m}}{d}{M}
\end{equation}
and use the surjectivity on the socles to obtain the result.
\end{proof}

\begin{theorem}\label{MainTheorem2}
Suppose $(A,\f{m})$ is a Noetherian local ring and let $M$ be a
finitely-generated $A$-module of dimension 2 such that $M$ has
positive depth, and \locoho{\f{m}}{1}{M} is finitely generated.
Then there exists an integer $\ell$ such that for every parameter
ideal \f{q} of $M$, if $\f{q} \sse \f{m}^{\ell}$ then the index of
reducibility of \f{q} on $M$ is given by
\begin{equation}
\index{A}{\f{q}}{M} = 2 \cdot \socdim{\locoho{\f{m}}{1}{M}} +
\socdim{\locoho{\f{m}}{2}{M}}.
\end{equation}
In particular, parameter ideals for $M$ in large powers of \f{m}
all have the same index of reducibility.
\end{theorem}

\begin{proof}
We begin by obtaining an integer $\ell$ from Proposition
\ref{T:lower bound}.  As in the proof of that proposition, we may
assume that $\ell$ is large enough so that every system of
parameters of $M$ contained in $\f{m}^{\ell}$ is a standard system
of parameters of $M$.

Let $a$, $b$ be a system of parameters of $M$ contained in
$\f{m}^{\ell}$ and set $\f{q} = (a, b)A$.  Then we have
\begin{equation}
\index{A}{\f{q}}{M} = \socdim{\frac{U_a + b M}{\f{q}M} + \frac{U_b
+ a M}{\f{q}M}} + \socdim{\locoho{\f{m}}{2}{M}},
\end{equation}
where $U_a = \unmixed{aA}{M} = \colon{aM}{M}{b}$, and similarly
for $U_b$.

Since $M$ has positive depth, $a$ and $b$ are regular elements on
$M$.  According to \cite[Theorem and Definition 17, p.
261]{Stuckrad.Vogel}, $a$ and $b$ both kill \locoho{\f{m}}{1}{M}.
Thus from the long exact sequence for local cohomology obtained
from the short exact sequence
\begin{equation}
\xymatrix{ 0 \ar[r] & M \ar[r]^{a} & M \ar[r] & M/aM \ar[r] & 0
 \\
}
\end{equation}
we see that $\locoho{\f{m}}{0}{M/aM} \cong \locoho{\f{m}}{1}{M}$.
Since $M/aM$ has dimension 1, \locoho{\f{m}}{0}{M/aM} is just
$\unmixed{aA}{M}/aM = \colon{aM}{M}{b}/aM$.  Furthermore, we have
a surjective homomorphism $\colon{aM}{M}{b} \to (\colon{aM}{M}{b}
+ bM)/bM$ whose kernel is $\colon{aM}{M}{b} \cap bM$.  Using the
definition of a standard system of parameters, we see that this
last expression is just $aM$.  Thus we have seen that
$\locoho{\f{m}}{1}{M} \cong (\colon{aM}{M}{b} + bM)/bM$.  The same
holds if we interchange $a$ and $b$.

At this point all that remains is to see that the sum
\begin{equation}
\frac{U_a + b M}{\f{q}M} + \frac{U_b + a M}{\f{q}M}
\end{equation}
is direct.  To this end, we note that
\begin{align}
(U_a + bM) \cap (U_b + aM) &= (U_a \cap (U_b + aM)) + bM\\
&= (U_a \cap U_b) + aM + bM\notag
\end{align}
Using the fact that $a$ and $b$ are regular on $M$ and form a
standard system of parameters of $M$, we see that
\begin{equation}
U_a \cap U_b = \colon{aM}{M}{b} \cap \colon{bM}{M}{a} = aM \cap
bM.
\end{equation}
Thus
\begin{equation}
(U_a + bM) \cap (U_b + aM) = \f{q}M,
\end{equation}
the sum is direct, and our proof is complete.
\end{proof}

\begin{corollary}
Let $(A, \f{m})$ be a Noetherian local ring of dimension at most
2.  Then $A$ is Gorenstein if and only if $A$ has finite local
cohomologies and every power of the maximal ideal contains an
irreducible parameter ideal.
\end{corollary}

\begin{proof}
If $A$ is Gorenstein, then all the local cohomology modules other
than \locoho{\f{m}}{2}{A} are zero, and all the parameter ideals
are irreducible.

For the other direction, note that it suffices to show that $A$ is
Cohen-Macaulay, since a Cohen-Macaulay local ring with an
irreducible parameter ideal is precisely a Gorenstein local ring.
The theorem is trivial in dimension zero: a zero dimensional
Noetherian local ring is Gorenstein if and only if the zero ideal
is irreducible. When the dimension of $A$ is one, the result
essentially goes back to Northcott-Rees \cite[Lemma 7]{NR-princ}:
If $A$ has SCI, then $A$ has positive depth.

In dimension 2, using the fact that the depth is positive, we know
from Theorem \ref{T:lower bound} that all parameter ideals in a
high power of the maximal ideal have the same index of
reducibility, namely $2 \cdot \socdim{\locoho{\f{m}}{1}{A}} +
\socdim{\locoho{\f{m}}{2}{A}}$. According to our hypothesis, this
integer must be 1.  Since \locoho{\f{m}}{2}{A} is a nonzero
Artinian module, it has a nonzero socle.  Thus is has a cyclic
socle, and \locoho{\f{m}}{1}{A} is an Artinian module with zero
socle.  Thus \locoho{\f{m}}{1}{A} is zero, so that $A$ is
Cohen-Macaulay.
\end{proof}

\begin{remark}
The technique in this section can be used to give a different
proof in the case of dimension 1.
\end{remark}

\section{An Example}

We begin with a lemma concerning idealizations.

\begin{definition}
Given a ring $R$ and an $R$-module $M$ we define the ring
\idealize{R}{M} to be the symmetric algebra of $M$ modulo the
square of the positive piece.  Thus \idealize{R}{M} is a graded
ring whose degree 0 piece is $R$, whose degree $1$ piece is $M$,
and whose components of degree greater than 1 are zero. This ring
is called the \emph{idealization} of $R$ with $M$, or the
\emph{trivial extension} of $R$ by $M$.
\end{definition}

For more information on idealization, we refer the reader to
\cite[p. 18]{Nagata} or \cite[Exercise 3.3.22]{Bruns.Herzog}.

\begin{lemma}\label{lemma}
\begin{enumerate}
\item Let $R$ be a ring and let $M$ be an $R$-module.  If $I$ is
an ideal of $R$, then
\[
\frac{\idealize{R}{M}}{I(\idealize{R}{M})} \cong
\idealize{\frac{R}{I}}{\frac{M}{IM}}.
\]

\item Let $(R,\f{m},k)$ be a Noetherian local ring and let $M$ be
an $R$-module.  Set $A = \idealize{R}{M}$.  Then
\[
\soc{A} = \soc{R} \cap \ann{R}{M} + \soc{M}.
\]
In particular, we have
\[
\socdim{A} = \socdim{M} + \Dim_k (\soc{R} \cap \ann{R}{M}).
\]

\item Let $(R,\f{m},k)$ be a Noetherian local ring and let $M$ be
an $R$-module.  Set $A = \idealize{R}{M}$.  If \f{q} is an
irreducible \f{m}-primary ideal of $R$ which does not contain
\ann{R}{M}, then $\f{q}A$ is a parameter ideal of $A$ for which
\begin{equation*}
\index{A}{\f{q}A}{A} = \index{R}{\f{q}}{M}+1.
\end{equation*}
\end{enumerate}
\end{lemma}

\begin{proof}
\begin{enumerate}
\item As $R$-modules, we have the isomorphism, and we immediately
see that this homomorphism respects multiplication.

\item We regard $A$ as $A = R+Mt$ with $t$ an indeterminate and
$t^2 = 0$. The maximal ideal of $A$ is $\f{m} + Mt$.  An element
$r + mt$ is in \soc{A} if and only if
\[
0 = (r + mt)(\f{m} + Mt) = r\f{m} + (rMt + \f{m}mt),
\]
which happens precisely when $r\in \soc{R} \cap \ann{R}{M}$ and $m
\in \soc{M}$, as desired.  The statement about dimensions now
follows from the fact that the $k$-vector space structure of the
socle is induced through its graded $R$-module structure.

\item We have
\begin{eqnarray*}
\index{A}{\f{q}A}{A} &=& \socdim{A/\f{q}A} \\ &=&
\socdim{\idealize{\frac{R}{\f{q}}}{\frac{M}{\f{q}M}}}
\\ &=&
\socdim{M/\f{q}M} + \Dim_k \soc{R/\f{q}} \cap
\ann{R/\f{q}}{M/\f{q}M}
\\ &=&
\index{R}{\f{q}}{M} + \Dim_k \soc{R/\f{q}}\cap
\ann{R/\f{q}}{M/\f{q}M}
\end{eqnarray*}

Now since \f{q} is an irreducible \f{m}-primary ideal, the socle
of $R/\f{q}$ is simple and essential, and is thus contained in
every nonzero submodule of $R/\f{q}$. Hence all that remains is to
note that \ann{R/\f{q}}{M/\f{q}M} is not zero, since it contains
$(\f{q} + \ann{R}{M})/\f{q}$.
\end{enumerate}
\end{proof}

We thank the referee for suggesting the following example in
dimension 2, and for pointing out the interesting paper of C. Lech
\cite{Lech}.

\begin{example}\label{example}
Let $k$ be a field and let $d > 1$ be an integer. Put $R = k[[x,
y, z_3, \ldots, z_d]]$, the formal power series ring over $k$ in
$d$ variables. Let $M = R/x^2 R$ and set $A = \idealize{R}{M}$.
 Let $w$ be a new variable; then $A \cong T/(x^2 w,w^2)T$, where $T
= k[[x,y,z_3, \ldots, z_d, w]]$. Thus $A$ is a complete Noetherian
local ring of dimension $d$ and depth $d - 1$.

For each integer $n \geq 3$ we define two parameter ideals of $R$:
\begin{eqnarray*}
\f{q} = (x^n, y^n, z_{3}^n, \ldots, z_{d}^n)R &\text{ and }
&\f{q}' = ((x + y)^n, xy^{n-1}, z_{3}^n, \ldots, z_{d}^n)R.
\end{eqnarray*}
Let $Q = \f{q}A$ and $Q' = \f{q}'A$; these are parameter ideals of
$A$ contained in $\f{m}^n$. We show that $\index{A}{Q}{A} = 2$
while $\index{A}{Q'}{A} = 3$.

Neither $Q$ nor $Q'$ contain $\ann{R}{M} = x^2 R$, and since $R$
is regular (hence Gorenstein), each of \f{q} and $\f{q}'$ are
irreducible. Hence from Part 3 of Lemma~\ref{lemma} we see that to
decide the index of reducibility of $Q$ and $Q'$, we just need to
consider
\[
\index{R}{\f{q}}{M} = \Dim_k R/(\f{q} + x^2 R),
\]
and similarly for $\f{q}'$.

We have
\begin{eqnarray*}
\f{q} + x^2 R = (x^2, y^n, z_{3}^n, \ldots, z_{d}^n)R &\text{ and
} &\f{q}' + x^2 R=(x^2, xy^{n-1}, y^n, z_{3}^n, \ldots, z_{d}^n)R.
\end{eqnarray*}
The first of these ideals is a parameter ideal for the Gorenstein
ring $R$, and is thus irreducible.  For the second we have
\[
(x^2, xy^{n - 1}, y^n, z_{3}^n, \ldots, z_{d}^n)R = (x^2, y^{n -
1}, z_{3}^n, \ldots, z_{d}^n)R \cap (x, y^n, z_{3}^n, \ldots,
z_{d}^n)R.
\]
Hence $\index{A}{Q}{A} = 1 + \index{R}{\f{q}}{M} = 2$ while
$\index{A}{Q'}{A} = 1 + \index{R}{\f{q}'}{M} = 3$.
\end{example}

\begin{remark}
We would like to point out that a bad example such as this may be
obtained among local Noetherian domains. The referee kindly
brought our attention to a paper of C. Lech \cite{Lech} which
proves that a complete Noetherian local ring $A$ is the completion
of a Noetherian local domain if and only if the prime ring of $A$
is a domain that acts without torsion on $A$, and the maximal
ideal of $A$ is not associated to $0$.  (Here the prime ring of a
ring is the subring generated by the multiplicative identity.)
\end{remark}


\begin{thebibliography}{99}

\bibitem[BH]{Bruns.Herzog}
W. Bruns and J. Herzog, \emph{Cohen-Macaulay Rings}, Cambridge
Studies in Advanced Mathematics, No. 39, Cambridge, 1993.

\bibitem[CHP]{Corso.Huneke.Vasconcelos} A. Corso, C. Huneke, and W. V. Vasconcelos, \emph{On the
integral closure of ideals}, manuscripta math. \textbf{95} (1998),
331-347.

\bibitem[CP]{Corso.Polini} A. Corso and C. Polini, \emph{Links of prime ideals and their
Rees algebras}, J. Alg. \textbf{178} (1995), 224-238.

\bibitem[CPV]{Corso.Polini.Vasconcelos} A. Corso, C. Polini, and W. V. Vasconcelos, \emph{Links of
prime ideals}, Math. Proc. Camb. Phil. Soc. \textbf{115} (1994),
431-436.

\bibitem[E]{Eisenbud}
D. Eisenbud, \emph{Commutative Algebra with a View Toward
Algebraic Geometry}, Springer-Verlag, New York, 1996.

\bibitem[EN]{Endo.Narita}
S. Endo and M. Narita, \emph{The Number of Irreducible Components
of an Ideal and the Semi-Regularity of a Local Ring}, Proc. Japan
Acad., \textbf{40} (1964), 627-630.

\bibitem[GSu]{Goto.Suzuki}
S. Goto and N. Suzuki, \emph{Index of Reducibility of Parameter
Ideals in a Local Ring}, Journal of Algebra, \textbf{87} (1984)
53-88.

\bibitem[GS1]{Goto.SakuraiI}
S. Goto and H. Sakurai, \emph{The equality $I^2 = QI$ in Buchsbaum
rings}, to appear in Rend. Sem. Mat. Univ. Padova \textbf{110}
(2003).

\bibitem[GS2]{Goto.SakuraiII}
S. Goto and H. Sakurai, \emph{The reduction exponent of socle
ideals associated to parameter ideals in a Buchsbaum local ring of
multiplicity two}, J. Math. Soc. Japan (to appear)

\bibitem[H]{Hochster} M. Hochster, \emph{Cyclic purity versus purity in
excellent Noetherian rings}, Trans. Amer. Math. Soc. \textbf{231}
(1977), 463-488.

\bibitem[K]{Kawasaki}
T. Kawasaki, \emph{On the Index of Reducibility of Parameter
Ideals and Cohen-Macaulayness in a Local Ring}, J. Math. Kyoto
Univ., \textbf{34-1} (1994), 219-226.

\bibitem[L]{Lech}
C. Lech, \emph{A method for constructing bad Noetherian local
rings}, \emph{Algebra, algebraic topology and their interactions}
(Stockholm, 1983), Lecture Notes in Math. \textbf{1183} (1986),
Springer-Verlag, Berlin $\cdot$ Heidelberg $\cdot$ New York
$\cdot$ Tokyo, 241-247

\bibitem[M]{Matsumura}
H. Matsumura, {\em Commutative Ring Theory}, Cambridge Studies in
Advanced Mathematics, No. 8, Cambridge, 1986.

\bibitem[Na]{Nagata}
M. Nagata, \emph{Local Rings}, Interscience Publishers, New York,
1962.

\bibitem[N]{Northc}
D. G. Northcott, {\em On Irreducible Ideals in Local Rings}, J.
London Math. Soc., \textbf{32} (1957), 82-88.

\bibitem[NR1]{NR-princ}
D. G. Northcott and D. Rees, \emph{Principal Systems}, Quart. J.
Math., \textbf{8} (1957), 119-27.

\bibitem[NR2]{NR-red}
D. G. Northcott and D. Rees, {\em Reductions of Ideals in Local
Rings,} Proc. Cambridge Phil. Soc., \textbf{50} (1954), 145-158.

\bibitem[ShV]{Sharpe.Vamos}
D. W. Sharpe and P. V. Vamos,  {\em Injective Modules}, Cambridge
University Press, Cambridge, 1972.

\bibitem[SV]{Stuckrad.Vogel}
J. St\"{u}ckrad and W. Vogel, \emph{Buchsbaum Rings and
Applications}, Springer-Verlag, Berlin, 1986

\end{thebibliography}
\end{document}